\documentclass{article}

\newtheorem{theorem}{Theorem}[section]
\newtheorem{corollary}[theorem]{Corollary}

\newtheorem{lemma}[theorem]{Lemma}
\newtheorem{proposition}[theorem]{Proposition}
\newtheorem{definition}[theorem]{Definition}

\evensidemargin\oddsidemargin
\begin{document}

\author{Vadim E. Levit and Eugen Mandrescu \\
%EndAName
Department of Computer Science\\
Holon Academic Institute of Technology\\
52 Golomb Str., P.O. Box 305\\
Holon 58102, ISRAEL\\
\{levitv, eugen\_m\}@barley.cteh.ac.il}
\title{A New Greedoid: \\
The Family of Local Maximum Stable Sets of a Forest}
\date{}
\maketitle

\begin{abstract}
A \textit{maximum stable set }in a graph $G$ is a stable set of maximum
cardinality. $S$ is a \textit{local maximum stable set} if it is a maximum
stable set of the subgraph of $G$ spanned by $S\cup N(S)$, where $N(S)$ is
the neighborhood of $S$. One theorem of Nemhauser and Trotter Jr. \cite
{NemhTro}, working as a useful sufficient local optimality condition for the
weighted maximum stable set problem, ensures that any local maximum stable
set of $G$ can be enlarged to a maximum stable set of $G$. In this paper we
demonstrate that an inverse assertion is true for forests. Namely, we show
that for any non-empty local maximum stable set $S$ of a forest $T$ there
exists a local maximum stable set $S_{1}$ of $T$, such that $S_{1}\subset S$
and $\left| S_{1}\right| =$ $\left| S\right| -1$. Moreover, as a further
strengthening of both the theorem of Nemhauser and Trotter Jr. and its
inverse, we prove that the family of all local maximum stable sets of a
forest forms a greedoid on its vertex set.
\end{abstract}

\section{Introduction}

Throughout this paper $G=(V,E)$ is a simple (i.e., a finite, undirected,
loopless and without multiple edges) graph with vertex set $V=V(G)$ and edge
set $E=E(G).$ If $X\subset V$, then $G[X]$ is the subgraph of $G$ spanned by 
$X$. By $G-W$ we mean the subgraph $G[V-W]$ , if $W\subset V(G)$. We also
denote by $G-F$ the partial subgraph of $G$ obtained by deleting the edges
of $F$, for $F\subset E(G)$, and we use $G-e$, if $W$ $=\{e\}$. The \textit{%
neighborhood} of a vertex $v\in V$ is the set $N(v)=\{w:w\in V$ \ \textit{and%
} $vw\in E\}$, whose cardinality is denoted by $deg(v)$. For $A\subset V$,
we denote $N(A,G)=\{v\in V-A:N(v)\cap A\neq \emptyset \}$ and $N[A,G]=A\cup
N(A)$, or shortly, $N(A)$ and $N[A]$, if no ambiguity. If $\left|
N(v)\right| =1$, then $v$ is a \textit{pendant vertex} of $G$. A stable set
of maximum size will be referred to as a \textit{maximum stable set} of $G$,
and the \textit{stability number }of $G$, denoted by $\alpha (G)$, is the
cardinality of a maximum stable set in $G$.

We call $A\subseteq V(G)$ a \textit{local maximum stable set} of $G$ if $%
A\in \Omega (G[A\cup N(A)])$. Let $\Omega (G)$ stand for the set $\{S:S$ 
\textit{is a maximum stable set of} $G\}$, and $\Psi (G)$ stand for the set
of all local maximum stable sets of graph $G$. For instance, any $A\subseteq
pend(G)$ is a local maximum stable set of $G$, where by $pend(G)$ we denote
the set of all pendant vertices of $G$. A graph $G$ is called $\alpha ^{+}$-%
\textit{stable} if $\alpha (G+e)=\nolinebreak \alpha (G)$, for any edge $%
e\in E(\overline{G})$, where $\overline{G}$ is the complement of $G$, \cite
{GunHarRall}. A \textit{matching} of $G$ is a set of edges no two of which
have a vertex in common. The \textit{matching number} $\mu (G)$ of $G$ is
the maximum size of a matching of $G$. A matching is \textit{perfect} if its
edges match up all vertices.

By $K_{n}$, $C_{n}$, $P_{n}$ we denote respectively, the complete graph on $%
n\geq 1$ vertices, the chordless cycle on $n\geq 4$ vertices, and the
chordless path on $n\geq 3$ vertices. Through all this paper we define a 
\textit{forest} as an acyclic graph of order greater than $1$, and a \textit{%
tree} as an acyclic connected graph of order greater than $1$. Since any
tree $T$ is also a bipartite graph, a well-known theorem of K\"{o}nig and
Egerv\'{a}ry assures that $\alpha (T)+\mu (T)=\left| V(T)\right| $, \cite
{Berge}, \cite{Eger}, \cite{Koen}. A \textit{perfect tree} is a tree having
a perfect matching, \cite{FrHedJaTre}. Gunther et al. proved in \cite
{GunHarRall}, that the perfect trees coincide with the $\alpha ^{+}$-stable
trees, and give also the following constructive characterization of $\alpha
^{+}$-stable trees:

\begin{theorem}
\label{th7}\cite{GunHarRall} $K_{2}$ is an $\alpha ^{+}$-stable tree. If $T$
is an $\alpha ^{+}$-stable tree, then the graph formed from $T$ by joining
one vertex of a new $K_{2}$ to some vertex of $T$ is also an $\alpha ^{+}$%
-stable tree.
\end{theorem}

In \cite{Zito} Zito extended some results of \cite{GunHarRall} and revealed
an elegant structure of maximum stable sets of a tree in terms of $\alpha $%
-critical edges, where an edge of a graph $G$ is called $\alpha $-critical
if $\alpha (G-e)>\alpha (G)$.

The following theorem concerning maximum stable sets in general graphs, due
to Nemhauser and Trotter Jr. \cite{NemhTro}, shows that for a special
subgraph $H$ of a graph $G$, some maximum stable set of $H$ can be enlarged
to a maximum stable set of $G$.

\begin{theorem}
\cite{NemhTro}\label{th3} Any local maximum stable set of a graph is a
subset of a maximum stable set.
\end{theorem}

Nemhauser and Trotter Jr. interpret this assertion as a sufficient local
optimality condition for a binary integer programming formulation of the
weighted maximum stable set problem, and use it to prove an impressive
result claiming that integer parts of solutions of the corresponding linear
programming relaxation retain the same values in the optimal solutions of
its binary integer programming counterpart. In other words, it means that a
well-known branch-and-bound heuristic for general integer programming
problems turns out to be an exact algorithm solving the weighted maximum
stable set problem.

Let us formulate an inverse version of Theorem \ref{th3} as follows:

\begin{description}
\item[Claim \{$k$\}.]  Any maximum stable set of a graph contains a local
maximum stable set of cardinality $k$.
\end{description}

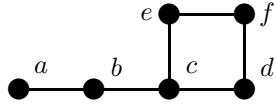
\begin{figure}[h]
\setlength{\unitlength}{1.0cm} 
\begin{picture}(5,1.5)\thicklines
  \multiput(4.5,0)(1,0){4}{\circle*{0.29}}
  \multiput(6.5,1)(1,0){2}{\circle*{0.29}}
  \put(4.5,0){\line(1,0){3}}
  \put(6.5,1){\line(1,0){1}}
  \multiput(6.5,0)(1,0){2}{\line(0,1){1}}
  \put(4.8,0.3){\makebox(0,0){$a$}}
  \put(5.8,0.3){\makebox(0,0){$b$}}
  \put(6.8,0.3){\makebox(0,0){$c$}}
  \put(7.8,0.3){\makebox(0,0){$d$}}
  \put(6.2,1){\makebox(0,0){$e$}}
  \put(7.8,1){\makebox(0,0){$f$}}
\end{picture}
\caption{A graph {with diverse local maximum stable sets}.}
\label{fig10}
\end{figure}

This claim is not valid for general graphs. For instance, Claim \{$k$\} is
false for all $k,$ $1\leq k<\alpha (G)$, if $G=C_{n},n\geq 4$. The graph $G$
in Figure \ref{fig10} shows another counterexample: any $S\in \Omega (G)$
contains some local maximum stable set, but these local maximum stable sets
are of different cardinalities. As examples, $\{a,c,f\}\in \Omega (G)$ but
only $\{a\}\in \Psi (G)$, while for $\{b,d,e\}\in \Omega (G)$ only $%
\{d,e\}\in \Psi (G)$.

Levit and Mandrescu proved in \cite{LevitMan} that any maximum stable set of
a tree $T$\ contains at least one of its pendant vertices (For any tree
Claim \{$k$\} is true for $k=1$), and if, in addition, $\alpha (T)>\left|
V(T)\right| /2$, then there exist at least two pendant vertices belonging to
all its maximum stable sets, i.e., in other words, any maximum stable set
includes both a local maximum stable set of size $1$ and size $2$ consisting
of pendant vertices (If the stability number of a tree is greater than half
of its order, then Claim \{$k$\} is true for $k=1,2$).

In this paper we prove that Claim \{$k$\} is true for any $k\in
\{1,2,...,\alpha (T)\}$, whenever $T$ is a forest. Moreover, we demonstrate
that for any $S\in \Omega (T)$, there is a chain 
\[
S_{1}\subset S_{2}\subset ...\subset S_{\alpha -1}\subset S_{\alpha }=S, 
\]
such that for all $1\leq i\leq \alpha (T)$, $\left| S_{i}\right| =i$ and $%
S_{i}$ is a local maximum stable set in $T$.

\begin{figure}[h]
\setlength{\unitlength}{1.0cm} 
\begin{picture}(5,1.5)\thicklines
  \multiput(4,0)(1,0){6}{\circle*{0.29}}
  \put(5,1){\circle*{0.29}}
  \multiput(7,1)(1,0){2}{\circle*{0.29}}
  \put(4,0){\line(1,0){5}}
  \put(7,1){\line(1,0){1}}
  \put(4,0){\line(1,1){1}}
  \put(5,0){\line(0,1){1}}
  \put(6,0){\line(-1,1){1}}
  \multiput(7,0)(1,0){2}{\line(0,1){1}}
  \put(9,0){\line(-1,1){1}}
\put(4,0.3){\makebox(0,0){$u$}}
\put(6,0.3){\makebox(0,0){$v$}}
\put(6.7,1){\makebox(0,0){$x$}}
\put(7.7,0.3){\makebox(0,0){$y$}}
\put(9,0.3){\makebox(0,0){$z$}}
\end{picture}
\caption{A graph with chains, which is not a forest.}
\label{fig11}
\end{figure}
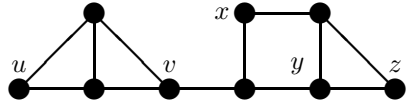

Notice that this property is not characteristic for forests only. The graph $%
G$ in Figure \ref{fig11} enjoys the same property, but it is not a forest.
Namely, $G$ has only two maximum stable sets, and each one of them generates
its corresponding chain: $\{u\}\subset \{u,v\}\subset \{u,v,z\}\subset
\{u,v,z,x\}$ and $\{u\}\subset \{u,v\}\subset \{u,v,y\}\subset \{u,v,y,x\}$.

In this form Claim \{$k$\} resembles an accessibility property of greedoids.
It turns out that this resemblance is not coincidental. Namely, we will
prove the following theorem.

\begin{theorem}
The family of local maximum stable sets of a forest forms a greedoid on its
vertex set.
\end{theorem}

The definition of greedoids we use in this paper is as follows.

\begin{definition}
\cite{BjZiegler}, \cite{KorLovSch} A greedoid is a pair $(E,\mathcal{F})$,
where $\mathcal{F}\subseteq 2^{E}$ is a set system satisfying the following
conditions:

\setlength {\parindent}{0.0cm}(\textit{Accessibility}) for every non-empty $%
X\in \mathcal{F}$ there is an $x\in X$ such that $X-\{x\}\in \mathcal{F}$;%
\setlength
{\parindent}{3.45ex}

\setlength {\parindent}{0.0cm}(\textit{Exchange}) for $X,Y\in \mathcal{F}%
,\left| X\right| =\left| Y\right| +1$, there is an $x\in X-Y$ such that $%
Y\cup \{x\}\in \mathcal{F}$.\setlength
{\parindent}{3.45ex}
\end{definition}

\section{The accessibility property}

\begin{lemma}
\label{lem2}If $A,B$ are two disjoint local maximum stable sets in $G$, such
that $A\cup B$ is stable, then $A\cup B$ is also a local maximum stable set
in $G$.
\end{lemma}

\setlength {\parindent}{0.0cm}\textbf{Proof.} Clearly, $A\cup B$ is stable
in $H=G[N[A]\cup N[B]]$.\setlength
{\parindent}{3.45ex}

Let $S\in \Omega (H)$. Hence, $\left| S\cap N[A]\right| \leq \left| A\right| 
$ and $\left| S\cap N[B]\right| \leq \left| B\right| $. Consequently, $%
\left| S\right| \leq \left| A\right| +\left| B\right| =\left| A\cup B\right|
\leq \left| S\right| $, and this implies that $A\cup B$ is a local maximum
stable set in $G$. \rule{2mm}{2mm}

\begin{theorem}
\label{th9}If $T$ is a perfect tree, then for any $S\in \Omega (T)$, there
is a chain 
\[
S_{1}\subset S_{2}\subset ...\subset S_{\alpha -1}\subset S_{\alpha }=S,
\]
such that for all $1\leq i\leq \alpha (T)$, $\left| S_{i}\right| =i$ and $%
S_{i}$ is a local maximum stable set in $T$.
\end{theorem}

\setlength {\parindent}{0.0cm}\textbf{Proof.} We use induction on $\alpha
(T) $. If $\alpha (T)=2$, then $T=P_{4}$, and the result is clear. Suppose
the assertion is true for perfect trees with stability number $\leq q$, and
let $T=(V,E)$ be a perfect tree with $\alpha (T)=q+1$. According to Theorem 
\ref{th7}, there is an edge $e=xy\in E$, such that $\left| N(y)\right|
=\left| \{x,w\}\right| =2$ and $x\in pend(T)$, because $T$ is $\alpha ^{+}$%
-stable, as well. Then, $T^{\prime }=T-\{x,y\}$ is also a perfect tree, and $%
\alpha (T^{\prime })=\alpha (T)-1=q$. If $S\in \Omega (T)$, then $S^{\prime
}=S-\{x,y\}\in \Omega (T^{\prime })$, and by induction hypothesis, there are 
\[
S_{1}\subset S_{2}\subset ...\subset S_{q-1}\subset S_{q}=S^{\prime }, 
\]
such that for all $1\leq i\leq \alpha (T^{\prime })$, $\left| S_{i}\right|
=i $ and $S_{i}$ is a local maximum stable set in $T^{\prime }$.%
\setlength
{\parindent}{3.45ex}

Suppose $S_{i}=\{v_{j}:1\leq j\leq i\},1\leq i\leq q$.

\textit{Case (}$\mathit{i}$\textit{).} $x\in S$. We show that 
\[
\{x\}\subset \{x\}\cup S_{1}\subset \{x\}\cup S_{2}\subset ...\subset
\{x\}\cup S_{q-1}\subset \{x\}\cup S_{q}=\{x\}\cup S^{\prime }=S 
\]
is a chain of local maximum stable sets of $T$, all included in $\{x\}\cup
S_{q}=\{x\}\cup S^{\prime }=S$ . If $y$ is not adjacent to any $v_{j},1\leq
j\leq q$, then the assertion is true by Lemma \ref{lem2} and the fact that
the neighborhoods of all $S_{i}$ in $T$ and in $T^{\prime }$ coincide.
Assume that $yv_{i}\in E(T)$, for some $i\in \{1,...,q\}$, i.e., $w=v_{i}$.
Then $W=$ $S_{i}\cup \{x\}$ is still a local maximum stable set of $T$,
because $W$ is stable and $\alpha (T[N[W]])=\left| S_{i}\right| +1=\left|
W\right| $.

\textit{Case (}$\mathit{ii}$)\textit{.} $y\in S$. Then $w\notin S^{\prime }$
and $S_{1}\subset S_{2}\subset ...\subset S_{q-1}\subset S_{q}=S^{\prime
}\subset S_{q+1}=S$ and all $S_{i}$ are local maximum stable sets in $%
T^{\prime }$, because the neighborhoods of all $S_{i}$ in $T$ and in $%
T^{\prime }$ coincide.

Thus, in both cases there exists a chain $S_{1}\subset S_{2}\subset
...\subset S_{q}\subset S_{q+1}=S$, such that for all $1\leq i\leq
q+1=\alpha (T)$, $\left| S_{i}\right| =i$ and $S_{i}$ is a local maximum
stable set in $T$. \rule{2mm}{2mm}

\begin{lemma}
\label{lem5}If $T_{1}$ is a subtree of the tree $T_{2}$, $A\subset V(T_{1})$%
, and $A\in \Psi (T_{2})$, then $A\in \Psi (T_{1})$.
\end{lemma}

\setlength {\parindent}{0.0cm}\textbf{Proof.} Let $N_{i}(A),i=1,2$, denote
the neighborhoods of $A$ in $T_{1},T_{2}$, respectively. Since $A$ is a
maximum stable set in $T_{2}[N_{2}[A]]$ and $N_{1}(A)\subseteq N_{2}(A)$, it
follows that $A$ is also a maximum stable set in $T_{1}[N_{1}[A]]$. \rule%
{2mm}{2mm}\setlength
{\parindent}{3.45ex}

\begin{lemma}
\label{lem7}Any tree $T_{1}$ can be embedded into a perfect tree $T_{2}$,
such that their stability numbers are equal.
\end{lemma}

\setlength {\parindent}{0.0cm}\textbf{Proof.} Let $M$ be a maximum matching
in $T_{1}$, and $V(M)$ be the vertices of $T_{1}$ matched by $M$. If $M$ is
a perfect matching, then we can choose $T_{2}=T_{1}$. Otherwise, if $%
\{v_{i}:1\leq i\leq q\}=V(T_{1})-V(M)$, let us define a new tree $T_{2}$ as
follows: 
\[
V(T_{2})=V(T_{1})\cup \{w_{i}:1\leq i\leq q\},E(T_{2})=E(T_{1})\cup
\{v_{i}w_{i}:1\leq i\leq q\}. 
\]
Clearly, $T_{2}$ is a perfect tree, since $M\cup \{v_{i}w_{i}:1\leq i\leq
q\} $ is a perfect matching in $T_{2}$, and $\mu (T_{2})=\left| M\right|
+q=\mu (T_{1})+q$. Consequently, by K\"{o}nig-Egerv\'{a}ry Theorem we obtain 
$\alpha (T_{2})=\left| V(T_{2})\right| -\mu (T_{2})=\left| V(T_{1})\right|
+q-\mu (T_{1})-q=\alpha (T_{1}),$and this completes the proof. \rule%
{2mm}{2mm}\setlength
{\parindent}{3.45ex}

\begin{figure}[h]
\setlength{\unitlength}{1.0cm} 
\begin{picture}(5,2.5)\thicklines
  \multiput(2,1)(1,0){10}{\circle*{0.29}}
  \multiput(3,0)(2,0){2}{\circle*{0.29}}
  \multiput(7,0)(1,0){2}{\circle*{0.29}}
  \multiput(10,0)(1,0){2}{\circle*{0.29}}
  \multiput(2,2)(1,0){2}{\circle*{0.29}}
  \put(5,2){\circle*{0.29}}
  \multiput(7,2)(1,0){5}{\circle*{0.29}}
  \put(2,1){\line(1,0){4}}
  \put(7,1){\line(1,0){4}}
  \put(2,1){\line(0,1){1}}
  \put(3,0){\line(0,1){2}}
  \put(5,0){\line(0,1){2}}
  \multiput(9,1)(0,0.2){5}{\circle*{0.09}}
  \multiput(11,1)(0,0.2){5}{\circle*{0.09}}
  \multiput(10,0)(0.2,0){5}{\circle*{0.09}} 
  \multiput(7,0)(0.2,0){5}{\circle*{0.09}}
  \multiput(8,0)(2,0){2}{\line(0,1){1}}
  \multiput(7,1)(1,0){2}{\line(0,1){1}}
  \put(10,1){\line(0,1){1}}
  \put(7,0){\circle*{0.34}}
  \put(9,2){\circle*{0.34}}
  \put(11,2){\circle*{0.34}}
  \put(11,0){\circle*{0.34}}
\end{picture}
\caption{A tree and one of its embeddings into a perfect tree (see the above
lemma). }
\label{fig101}
\end{figure}
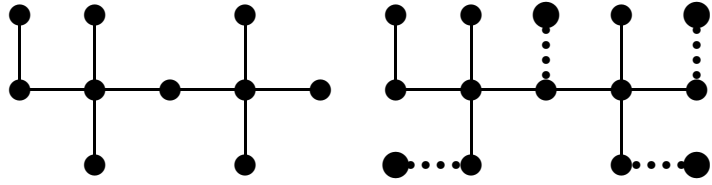

\begin{proposition}
\label{lem6}Any tree contains a maximum matching covering all its internal
vertices.
\end{proposition}

\setlength {\parindent}{0.0cm}\textbf{Proof.} Let $M$ be a maximum matching
in a tree $T$, and suppose that some vertex $v\in V(T)-pend(T)$ is not
matched. Hence, all neighbors of $v$ are matched, otherwise $M$ is not a
maximum matching. If $u\in N(v)$ and $uw\in M$, then $M_{1}=M\cup
\{vu\}-\{uw\}$ is also a maximum matching. If $w\in pend(T)$, we continue
with another internal vertex of $T$, unmatched by $M_{1}$, if such a vertex
exists. If $w\notin pend(T)$, then all its neighbors are matched by $M_{1}$,
and we can choose a vertex $x\in N(w)-\{u\}$, for which some edge $xy\in
M_{1}$. Hence $M_{2}=M_{1}\cup \{wx\}-\{xy\}$ is again a maximum matching in 
$T$. If $y\notin pend(T)$, we continue in the same manner, until some
pendant vertex stops us. The final matching $M_{p}$ saturates $v$ and all
the internal vertices matched by $M$. If there exists in $T$ an internal
vertex $a$ still unmatched by $M_{p}$, we repeat the procedure. After a
finite number of steps, we obtain a maximum matching covering all the
internal vertices of $T$. \rule{2mm}{2mm}\setlength
{\parindent}{3.45ex}\newline

Combining Proposition \ref{lem6} and Lemma \ref{lem7} we obtain the
following corollary.

\begin{corollary}
\label{cor1}Any tree $T_{1}$ can be embedded into a perfect tree $T_{2}$,
such that all the new edges are adjacent to pendant vertices of $T_{1}$, and 
$\alpha (T_{1})=\alpha (T_{2})$.
\end{corollary}

\begin{figure}[h]
\setlength{\unitlength}{1.0cm} 
\begin{picture}(5,2.5)\thicklines
  \multiput(2,1)(1,0){10}{\circle*{0.29}}
  \multiput(3,0)(2,0){2}{\circle*{0.29}}
  \multiput(7,0)(1,0){2}{\circle*{0.29}}
  \multiput(10,0)(1,0){2}{\circle*{0.29}}
  \multiput(2,2)(1,0){2}{\circle*{0.29}}
  \put(5,2){\circle*{0.29}}
  \multiput(7,2)(1,0){5}{\circle*{0.29}}
  \put(2,1){\line(1,0){4}}
  \put(7,1){\line(1,0){4}}
  \put(2,1){\line(0,1){1}}
  \put(3,0){\line(0,1){2}}
  \put(5,0){\line(0,1){2}}
  \multiput(7,0)(0.2,0){5}{\circle*{0.09}}
  \multiput(10,0)(0.2,0){5}{\circle*{0.09}} 
  \multiput(8,0)(2,0){2}{\line(0,1){1}}
  \multiput(7,1)(1,0){2}{\line(0,1){1}}
  \put(10,1){\line(0,1){1}}
  \multiput(9,2)(0.2,0){5}{\circle*{0.09}}
  \multiput(11,1)(0,0.2){5}{\circle*{0.09}} 
  \put(7,0){\circle*{0.34}}
  \put(9,2){\circle*{0.34}}
  \put(11,2){\circle*{0.34}}
  \put(11,0){\circle*{0.34}}
\end{picture}
\caption{A tree and one of its embeddings into a perfect tree (see the above
corollary).}
\label{fig102}
\end{figure}
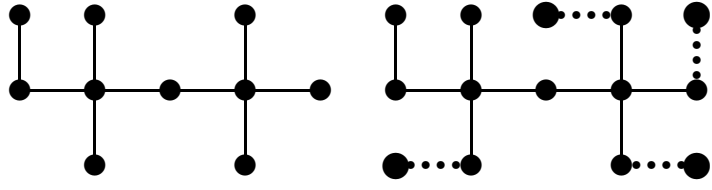

\begin{theorem}
\label{th10}If $T$ is a non-perfect tree, then for any $S\in \Omega (T)$,
there exists a chain $S_{1}\subset S_{2}\subset ...\subset S_{\alpha
-1}\subset S_{\alpha }=S,$such that for all $1\leq i\leq \alpha (T)$, $%
\left| S_{i}\right| =i$ and $S_{i}$ is a local maximum stable set in $T$.
\end{theorem}

\setlength {\parindent}{0.0cm}\textbf{Proof.} According to Lemma \ref{lem7}, 
$T$ can be embedded into a perfect tree $T^{\prime }$, such that $\alpha
(T)=\alpha (T^{\prime })$. If $S\in \Omega (T)$, it follows that $S\in
\Omega (T^{\prime })$, and by Theorem \ref{th9}, there is a chain $%
S_{1}\subset S_{2}\subset ...\subset S_{\alpha -1}\subset S_{\alpha }=S$,
such that for all $1\leq i\leq \alpha (T^{\prime })$, $\left| S_{i}\right|
=i $ and $S_{i}$ is a local maximum stable set in $T^{\prime }$. Hence,
Lemma \ref{lem5} ensures that all $S_{i},1\leq i\leq \alpha (T^{\prime
})=\alpha (T)$, are local maximum stable sets in $T$ as well. \rule{2mm}{2mm}%
\setlength
{\parindent}{3.45ex}

\begin{proposition}
\label{Chain}If $T$ is a tree, $S\in \Psi (T)$, and $\left| S\right| =k$,
then there exists a chain $S_{1}\subset S_{2}\subset ...\subset
S_{k-1}\subset S_{k}=S,$such that for all $1\leq i\leq k$, $\left|
S_{i}\right| =i$ and $S_{i}$ is a local maximum stable set in $T$.
\end{proposition}

\setlength {\parindent}{0.0cm}\textbf{Proof.} Suppose that $T_{1}=T[S\cup
N(S)]$ is also a tree. According to Theorems \ref{th9} and \ref{th10}, it
follows that there is a chain $S_{1}\subset S_{2}\subset ...\subset
S_{k-1}\subset S_{k}=S$, such that for all $1\leq i\leq k$, $\left|
S_{i}\right| =i$ and $S_{i}$ is a local maximum stable set in $T_{1}$. Since 
$N[S,T]=N[S,T_{1}]$, any $S_{i}$ has $N[S_{i},T_{1}]\subseteq N[S,T]$, and
therefore we get $N[S_{i},T_{1}]=N[S_{i},T]$. Consequently, all $S_{i},1\leq
i\leq k$, are local maximum stable sets in $T$. Assume $T_{1}$ is a forest.
Without loss of generality, we may suppose that $T_{1}$ contains two trees $%
T_{2}$ and $T_{3}$. Then $S_{2}=S\cap V(T_{2})$ and $S_{3}=S\cap V(T_{3})$
are also in $\Psi (T)$, and as above, there are two chains of local maximum
stable sets in $T$, as follows: \setlength
{\parindent}{3.45ex} 
\[
S_{1}^{^{\prime }}\subset S_{2}^{^{\prime }}\subset ...\subset
S_{k_{1}-1}^{^{\prime }}\subset S_{k_{1}}^{^{\prime }}=S_{2}\ and\
S_{1}^{^{\prime \prime }}\subset S_{2}^{^{\prime \prime }}\subset ...\subset
S_{k_{2}-1}^{^{\prime \prime }}\subset S_{k_{2}}^{^{\prime \prime }}=S_{3}. 
\]

Then using Lemma \ref{lem2}, we get a chain for $S$ itself, namely:

\[
S_{1}^{^{\prime }}\subset S_{2}^{^{\prime }}\subset ...\subset
S_{k_{1}}^{^{\prime }}\subset S_{k_{1}}^{^{\prime }}\cup S_{1}^{^{\prime
\prime }}\subset S_{k_{1}}^{^{\prime }}\cup S_{2}^{^{\prime \prime }}\subset
...\subset S_{k_{1}}^{^{\prime }}\cup S_{k_{2}}^{^{\prime \prime
}}=S_{2}\cup S_{3}=S_{3}, 
\]
and this completes the proof. \rule{2mm}{2mm}\newline

The following accessibility property for the family of local maximum stable
sets of a tree is an equivalent form of Proposition \ref{Chain}.

\begin{theorem}[Accessibility Property for Trees]
\label{Accessibility}If $S\in \Psi (T)$ and $T$ is a tree, then there exists
some $S_{1}\subset S$, such that $S_{1}\in \Psi (T)$ and $\left|
S_{1}\right| =\left| S\right| -1$.
\end{theorem}

Notice that if $T$ is a forest and $\{T_{i}:1\leq i\leq q\}$ are its
connected components, then $\Psi (T)=\cup \{\Psi (T_{i}):1\leq i\leq q\}$,
and using Lemma \ref{lem2} and Theorem \ref{Accessibility}, we obtain:

\begin{theorem}[Accessibility Property for Forests]
\label{Accessibilityf}If $S\in \Psi (T)$ and $T$ is a forest, then there
exists some $S_{1}\subset S$, such that $S_{1}\in \Psi (T)$ and $\left|
S_{1}\right| =\left| S\right| -1$.
\end{theorem}

\begin{figure}[h]
\setlength{\unitlength}{1.0cm} 
\begin{picture}(5,1.5)\thicklines
  \multiput(4,0)(1,0){6}{\circle*{0.29}}
  \multiput(4,1)(1,0){2}{\circle*{0.29}}
  \multiput(7,1)(1,0){3}{\circle*{0.29}}
  \put(4,0){\line(1,0){5}}
  \put(7,1){\line(1,0){1}}
  \multiput(4,0)(1,0){2}{\line(0,1){1}}
  \multiput(7,0)(2,0){2}{\line(0,1){1}}
  \put(5.3,1){\makebox(0,0){$a$}}  
  \put(4.25,0.3){\makebox(0,0){$b$}}
  \put(6.2,0.3){\makebox(0,0){$c$}}
  \put(6.7,1){\makebox(0,0){$d$}}
  \put(9.3,1){\makebox(0,0){$e$}}
\end{picture}
\caption{$\{a\}\subset \{a,b\}\subset \{a,b,c\}\subset \{a,b,c,d\}\subset
\{a,b,c,d,e\}$.}
\label{fig40}
\end{figure}
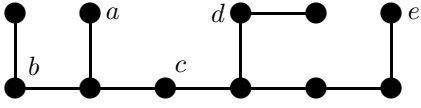

Figure \ref{fig40} presents a chain of local maximum stable sets in a tree.

\section{The exchange property}

According to Theorem \ref{th3} of Nemhauser and Trotter Jr. \cite{NemhTro},
any $A\in \Psi (G)$ can be enlarged to some $S\in \Omega (G)$. We show that
for every $S\in \Omega (G)$ this enlargement can be implemented using only
elements of $S$.

\begin{theorem}[Exchange Version of Nemahauser's and Trotter's Theorem]
\label{Exchange Version of N&T }If $S_{2}\in \Omega (G)$ and $S_{1}\in \Psi
(G)$, then there exists $S_{3}\subseteq S_{2}-S_{1}$ such that $S_{1}\cup
S_{3}\in \Omega (G)$.
\end{theorem}

\setlength {\parindent}{0.0cm}\textbf{Proof.} Since $S_{1}\in \Psi (G)$, it
follows that $\left| N[S_{1}]\cap S_{2}\right| \leq \left| S_{1}\right| $,
and consequently $S_{3}=S_{2}-S_{2}\cap N[S_{1}]$ is stable and $\left|
S_{3}\right| \geq \left| S_{2}\right| -\left| S_{1}\right| $. Hence we get
that $S_{1}\cup S_{3}$ is stable and $\alpha (G)\geq \left| S_{1}\cup
S_{3}\right| =\left| S_{1}\right| +\left| S_{3}\right| \geq \left|
S_{1}\right| +\left| S_{2}\right| -\left| S_{1}\right| =\alpha (G)$.
Therefore, $\alpha (G)=\left| S_{1}\cup S_{3}\right| $, i.e., $S_{1}\cup
S_{3}\in \Omega (G)$. \rule{2mm}{2mm}\setlength
{\parindent}{3.45ex}.

\begin{corollary}[$\left( \alpha -1,\alpha \right) $ Exchange Property]
\label{(Alpha-1,Alpha) Exchange Property}If $S_{2}\in \Omega (G),S_{1}\in
\Psi (G)$ and $\left| S_{2}\right| =\left| S_{1}\right| +1$, then there
exists $v\in S_{2}-S_{1}$ such that $S_{1}\cup \{v\}\in \Omega (G)$.
\end{corollary}

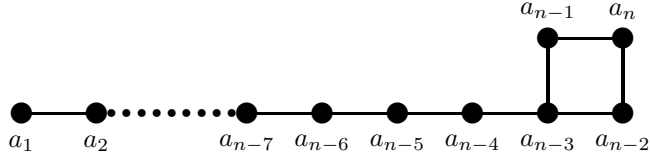
\begin{figure}[h]
\setlength{\unitlength}{1.0cm} 
\begin{picture}(5,2)\thicklines
  \multiput(2.5,0.5)(1,0){2}{\circle*{0.29}}
  \multiput(5.5,0.5)(1,0){6}{\circle*{0.29}}
  \multiput(9.5,1.5)(1,0){2}{\circle*{0.29}}
  \put(2.5,0.5){\line(1,0){1}}
  \put(5.5,0.5){\line(1,0){5}}
  \put(9.5,1.5){\line(1,0){1}}
  \multiput(3.5,0.5)(0.2,0){11}{\circle*{0.09}}
  \multiput(9.5,0.5)(1,0){2}{\line(0,1){1}}
  \put(2.5,0.1){\makebox(0,0){$a_{1}$}}
  \put(3.5,0.1){\makebox(0,0){$a_{2}$}}
  \put(5.5,0.1){\makebox(0,0){$a_{n-7}$}}
  \put(6.5,0.1){\makebox(0,0){$a_{n-6}$}}
  \put(7.5,0.1){\makebox(0,0){$a_{n-5}$}}
  \put(8.5,0.1){\makebox(0,0){$a_{n-4}$}}
  \put(9.5,0.1){\makebox(0,0){$a_{n-3}$}}
  \put(10.5,0.1){\makebox(0,0){$a_{n-2}$}}
  \put(9.5,1.85){\makebox(0,0){$a_{n-1}$}}
  \put(10.5,1.85){\makebox(0,0){$a_{n}$}}
\end{picture}
\caption{$\left( \alpha -2,\alpha -1\right) $ and $\left( 1,2\right) $
counterexamples to the exchange property.}
\label{fig70}
\end{figure}

Let us notice that if $S_{1},S_{2}\in \Psi \left( G\right) -\Omega (G)$,
then sometimes there is no $v\in S_{2}-S_{1}$ such that $S_{1}\cup \{v\}\in
\Psi (G)$. For instance, for the graph $G$ in Figure \ref{fig70} we have $%
\{a_{1}\},\{a_{n-2},a_{n-1}\}\in \Psi (G)$, but $\{a_{1},a_{n-2}\},%
\{a_{1},a_{n-1}\}\notin \Psi (G)$, provided $n\geq 6$. This example shows
that for every $n\geq 6$ there exists a graph $G$ of order $n$ with a pair
of local maximum stable sets of cardinalities $1$ and $2$ for which the
exchange property is not valid.

For the graph $G$ in Figure \ref{fig70}, if $n\geq 8$ is even then $2\alpha
\left( G\right) =n$. It is easy to check that 
\begin{eqnarray*}
S_{1} &=&\left\{ a_{1},a_{3},...,a_{n-5}\right\} ,S_{2}=\left\{
a_{1},a_{3},...,a_{n-7},a_{n-2},a_{n-1}\right\} \in \Psi (G), \\
\left| S_{1}\right| &=&\alpha \left( G\right) -2,\left| S_{2}\right| =\alpha
\left( G\right) -1,S_{1}\cup \left\{ a_{n-2}\right\} ,S_{1}\cup \left\{
a_{n-1}\right\} \notin \Psi (G).
\end{eqnarray*}
It means that for every even $n\geq 8$ there exists a graph $G$ of order $n$
with a pair of local maximum stable sets of cardinalities $\alpha \left(
G\right) -2$ and $\alpha \left( G\right) -1$ for which the exchange property
is not valid.

The next theorem shows that for forests the exchange property, i.e., the
assertion in Corollary \ref{(Alpha-1,Alpha) Exchange Property}, is true even
if the local maximum stable set $S_{2}\notin \Omega (G)$.

\begin{theorem}[Exchange Property]
\label{Exchange Property}If $S_{1},S_{2}\in \Psi (T)$ and $\left|
S_{2}\right| =\left| S_{1}\right| +1$, then there exists $v\in S_{2}-S_{1}$,
such that $S_{1}\cup \{v\}\in \Psi (T)$.
\end{theorem}

\setlength {\parindent}{0.0cm}\textbf{Proof.} We use induction on $k=\left|
S_{2}\right| $. \setlength
{\parindent}{3.45ex}

Base 1. If $\left| S_{2}\right| =\left| \{v\}\right| =1$, then $v\in pend(T)$%
, $S_{1}=\emptyset $, and clearly $S_{1}\cup \{v\}\in \Psi (T)$.

Base 2. If $\left| S_{2}\right| =\left| \{v,w\}\right| =2$, then at least
one of $v,w$, say $v$, is in $pend(T)$ (see Proposition \ref{Chain}). Let $%
\left| S_{1}\right| =\left| \{u\}\right| =1$. If $u\in S_{2}$, then either $%
S_{1}\cup \{v\}=S_{2}\in \Psi (T)$, or $S_{1}\cup \{w\}=S_{2}\in \Psi (T)$.
If $u\notin S_{2}$, then $S_{1}\cup \{v\}\in \Psi (T)$ according to Lemma 
\ref{lem2}.

Suppose that the assumption is true for sets of cardinality $\leq k$, and
let $S_{2}$ be of cardinality $k+1$. Since any local maximum stable set
contains at least one pendant vertex (according to Proposition \ref{Chain}),
we distinguish between the following three cases: ($\emph{i}$) $S_{2}$
contains some vertex $v\in pend(T)-N[S_{1}]$; ($\emph{ii}$) $S_{2}$ contains
some vertex $v\in pend(T)\cap N(S_{1})$; ($\emph{iii}$) there is some $v\in
pend(T)\cap S_{2}\cap S_{1}$.

\textit{Case} ($\emph{i}$) There exists some $v\in S_{2}\cap pend(T)-N\left[
S_{1}\right] $. Since $\{v\}\cup S_{1}$ is stable, Lemma \ref{lem2} implies
that $S_{1}\cup \{v\}\in \Psi (T)$.

\begin{figure}[h]
\setlength{\unitlength}{1.0cm} 
\begin{picture}(5,2.3)\thicklines
  \multiput(6,2)(1,0){2}{\circle*{0.29}}
  \multiput(4,1)(1,0){5}{\circle*{0.29}}
  \multiput(6,0)(1,0){2}{\circle*{0.29}}
  \put(6,2){\line(1,0){1}}
  \put(6,0){\line(1,0){1}}
  \put(4,1){\line(1,0){4}}
  \put(7,0){\line(0,1){2}}
  \put(4,1.35){\makebox(0,0){$a$}}
  \put(6,1.35){\makebox(0,0){$b$}}
  \put(8,1.35){\makebox(0,0){$e$}}
  \put(6,2.35){\makebox(0,0){$c$}}
  \put(7,2.35){\makebox(0,0){$d$}}
\end{picture}
\caption{Tree $T$ with $\{a,b,c,e\},\{a,d,e\}\in \Psi (T)$, and also $%
\{a,d,e\}\cup \{b\}\in \Psi (T)$.}
\label{fig43}
\end{figure}
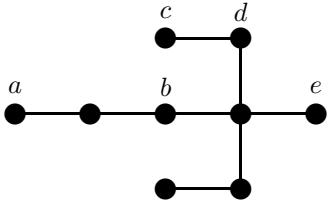

\textit{Case }($\emph{ii}$) There exists some $v\in S_{2}\cap pend(T)$, such
that $N(v)\cap S_{1}=\{u\}$. Figure \ref{fig43} illustrates this case.

Firstly, we show that $S_{2}-\{v\},S_{1}-\{u\}\in \Psi (T-\{u,v\})$. If $%
S_{2}-\{v\}$ does not belong to $\Psi (T-\{u,v\})$, then there is $A\in \Psi
(N[S_{2}-\{v\}])$ with $\left| A\right| >\left| S_{2}-\{v\}\right| $, and
therefore $A\cup \{v\}$ is a stable set in $N[S_{2}]$ larger than $S_{2}$,
in contradiction with $S_{2}\in \Psi (T)$. If $S_{1}-\{u\}\notin \Psi
(T-\{u,v\})$, then there exists $A\in \Psi (N[S_{1}-\{u\}])$ with $\left|
A\right| >\left| S_{1}-\{u\}\right| $, and therefore $A\cup \{u\}$ is a
stable set in $N[S_{1}]$ larger than $S_{1}$, in contradiction with $%
S_{1}\in \Psi (T)$. By the induction hypothesis, there exists $x\in
(S_{2}-\{v\})-(S_{1}-\{u\})$ such that $(S_{1}-\{u\})\cup \{x\}\in \Psi
(T-\{u,v\})$.

Secondly, to complete the proof of the theorem for the case ($\emph{ii}$) we
will show that $S_{1}\cup \{x\}\in \Psi (T)$. The set $\{u,x\}$ is stable,
since otherwise $(S_{1}-\{u\})\cup \{v,x\}$ is a stable set in $N[S_{1}]$
with its cardinality larger than the cardinality of $S_{1}$, in
contradiction with $S_{1}\in \Psi (T)$. Consequently, $S_{1}\cup \{x\}$ is
stable, because $(S_{1}-\{u\})\cup \{x\}$ is also stable. Since $%
(S_{1}-\{u\})\cup \{x\}$ is a maximum stable set in 
\[
N[(S_{1}-\{u\})\cup \{x\}]=N[S_{1}\cup \{x\}]-\{u,v\} 
\]
and $S_{1}\cup \{x\}$ is stable, it follows that $S_{1}\cup \{x\}$ is a
maximum stable set in $N[S_{1}\cup \{x\}]$, i.e., $S_{1}\cup \{x\}\in \Psi
(T)$.

\begin{figure}[h]
\setlength{\unitlength}{1.0cm} 
\begin{picture}(5,1)\thicklines
  \multiput(3,0)(1,0){8}{\circle*{0.29}}
  \put(3,0){\line(1,0){7}}
  \put(3,0.4){\makebox(0,0){$a$}}
  \put(4,0.4){\makebox(0,0){$b$}}
  \put(5,0.4){\makebox(0,0){$c$}}
  \put(6,0.4){\makebox(0,0){$d$}}
  \put(7,0.4){\makebox(0,0){$e$}}
  \put(8,0.4){\makebox(0,0){$f$}}
  \put(9,0.4){\makebox(0,0){$g$}}
  \put(10,0.4){\makebox(0,0){$h$}}
  \end{picture}
\caption{Tree $T$ with $\{a,c,e,h\},\{a,f,h\}\in \Psi (T)$, and also $%
\{a,f,h\}\cup \{c\}\in \Psi (T)$.}
\label{fig42}
\end{figure}
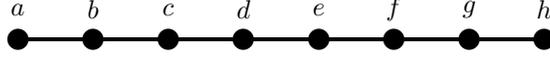

\textit{Case }($\emph{iii}$) There exists some $v\in pend(T)\cap S_{2}\cap
S_{1}$. Figure \ref{fig42} illustrates this case.

If $N(v)=\{u\}$, then $S_{2}-\{v\},S_{1}-\{v\}\in \Psi (T-\{u,v\})$. By the
induction hypothesis, there exists $x\in (S_{2}-\{v\})-(S_{1}-\{v\})$ such
that $(S_{1}-\{v\})\cup \{x\}\in \Psi (T-\{u,v\})$. We show now that $%
S_{1}\cup \{x\}\in \Psi (T)$. Clearly $x\neq u$, and therefore, $\{v,x\}$ is
stable. Consequently, $S_{1}\cup \{x\}$ is also stable, because $%
(S_{1}-\{v\})\cup \{x\}$ is stable. Since $(S_{1}-\{v\})\cup \{x\}$ is a
maximum stable set in 
\[
N[(S_{1}-\{v\})\cup \{x\}]=N[S_{1}\cup \{x\}]-\{v\} 
\]
and $S_{1}\cup \{x\}$ is stable, it follows that $S_{1}\cup \{x\}$ is a
maximum stable set in $N[S_{1}\cup \{x\}]$, i.e., $S_{1}\cup \{x\}\in \Psi
(T)$. \rule{2mm}{2mm}\newline

\begin{figure}[h]
\setlength{\unitlength}{1.0cm} 
\begin{picture}(5,1)\thicklines
  \multiput(4,0)(1,0){6}{\circle*{0.29}}
  \put(4,0){\line(1,0){5}}
  \put(4,0.4){\makebox(0,0){$a$}}
  \put(5,0.4){\makebox(0,0){$b$}}
  \put(6,0.4){\makebox(0,0){$c$}}
  \put(7,0.4){\makebox(0,0){$d$}}
  \put(8,0.4){\makebox(0,0){$e$}}
  \put(9,0.4){\makebox(0,0){$f$}}

  \end{picture}
\caption{Tree $T$ with $S_{2}=\{a,c\},S_{1}=\{f\}\in \Psi (T)$. For $a\in
S_{2}-S_{1},\{a,f\}\in \Psi (T)$; nevertheless, for $c\in S_{2}-S_{1}$, $%
\{c,f\}$ is stable but does not belong to $\Psi (T)$. }
\label{fig41}
\end{figure}
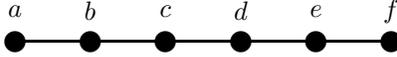

The example in Figure\ \ref{fig41} shows that even for trees there exist a
pair of local maximum stable sets $S_{1},S_{2}$, such that $S_{1}\cup \{v\}$
is not a local maximum stable set for all $v\in S_{2}-S_{1}$.

Combining Theorems \ref{Accessibilityf} and \ref{Exchange Property}, we
finally obtain:

\begin{theorem}
The family of local maximum stable sets of a forest forms a greedoid on its
vertex set.
\end{theorem}

\section{Conclusions}

In this paper we have proved that an inverse statement to Theorem \ref{th3}
due to Nemhauser and Trotter is true for forests. As a further strengthening
of both Theorem \ref{th3} and its inverse, we also have shown that the
family of local maximum stable sets of a forest generates a greedoid on its
vertex set. It seems to us quite interesting to find a general description
of such a greedoid. This also bring us to the following open problems: for
which classes of graphs an inverse of Theorem \ref{th3} is still true, and
for which classes of graphs their families of local maximum stable sets form
greedoids?

\section{Acknowledgment}

We would like to thank Endre Boros for drawing our attention to the paper of
Nemhauser and Trotter Jr., and for many helpful discussions.

\end{document}